\documentclass[pra, showpacs, showkeys, onecolumn, amsmath, amssymb, superscriptaddress, preprintnumbers, floatfix, preprint]{revtex4}

\usepackage{bm,amsfonts, mathtools}
\usepackage{graphicx,color, wasysym}
\usepackage{textcomp}
\usepackage{hyperref}

\begin{document}

\title{Definite integrals and operational methods}

\author{\large{D.~Babusci}}
\email{danilo.babusci@lnf.infn.it}

\affiliation{INFN - Laboratori Nazionali di Frascati, v.le E. Fermi, 40, I 00044 Frascati (Roma), Italy}

\author{\large{G.~Dattoli}}
\email{dattoli@frascati.enea.it}

\affiliation{ENEA - Centro Ricerche Frascati, v.le E. Fermi, 45, I 00044 Frascati (Roma), Italy}

\author{\large{G.~H.~E.~Duchamp}}
\email{ghed@lipn-univ.paris13.fr}

\affiliation{Universit\'e Paris XIII, LIPN, Institut Galil\'{e}e, CNRS UMR 7030, 99 Av. J.-B. Clement, F 93430 Villetaneuse, France} 

\author{\large{K.~G\'{o}rska}}
\email{kasia\_gorska@o2.pl}

\affiliation{Instituto de F\'{\i}sica, Universidade de S\~{a}o Paulo, P.O.Box 66318, BR 05315-970 S\~{a}o Paulo, SP, Brasil}

\affiliation{H. Niewodnicza\'{n}ski Institute of Nuclear Physics, Polish Academy of Sciences, ul.Eljasza-Radzikowskiego 152, PL 31342 Krak\'{o}w, Poland}

\author{\large{K.~A.~Penson}}
\email{penson@lptl.jussieu.fr}

\affiliation{Laboratoire de Physique Th\'{e}orique de la Mati\`{e}re Condens\'{e}e (LPTMC), Universit\'{e} Pierre et Marie Curie, CNRS UMR 7600, Tour 13 - 5i\`{e}me \'{e}t., Bo\^{i}te Courrier 121, 4 place Jussieu, F 75252 Paris Cedex 05, France}


\begin{abstract}
An operational method, already employed to formulate a generalization of the Ramanujan master theorem, is applied to the evaluation of integrals of various types. This technique provides a very flexible and powerful tool yielding new results encompassing different aspects of the special function theory.
\end{abstract}

\maketitle


The use of operational methods is, sometimes, very much useful to "guess" specific theorems in various fields of analysis. The drawback of such an approach is that any of its consequences should be validated by a more rigorous procedure. This is indeed the case of the Ramanujan Master Theorem \cite{BBerndt40, GHHardy40}, according to which if a function $f(x)$ can be expanded around $x = 0$ in a power series of the form
\begin{equation}\label{eq1}
f(x) = \sum_{n=0}^{\infty} \frac{\varphi(n)}{n!} (-x)^{n},
\end{equation}
with $\varphi(0) = f(0) \neq 0$, then
\begin{equation}\label{eq2}
\int_{0}^{\infty} x^{\nu-1} f(x)\, dx = \Gamma(\nu)\, \varphi(\nu),
\end{equation}
where $\Gamma(z)$ is the Euler's Gamma function. The rigorous proof of this identity has been obtained in Ref. \cite{TAmdeberhan09, TAmdeberhan11}. The simpler but heuristic proof is given in \cite{DBabusci11, KGorska11}. It is based on setting
\begin{equation}\label{eq3}
f(x) = \sum_{n=0}^{\infty} \frac{\hat{c}^{\,n}}{n!} (-x)^{n} \varphi(0) = e^{-\hat{c}\, x} \varphi(0),
\end{equation}
with $\hat{c}$ being an umbral operator such that \cite{SRoman05}
\begin{equation}\label{eq4}
\hat{c}^{\,n} \varphi(0) = \varphi(n).
\end{equation}
The theorem has been exploited in \cite{DBabusci11, KGorska11} in a generalized version to obtain, in a very simple way, old and new formulas for integrals and successive derivatives of Bessel functions. The noticeable feature emerging from such a procedure is that cylindrical Bessel functions can be formally treated as functions involving the operator $\hat{c}$. For example, the $n$-th order cylindrical Bessel function can be written as
\begin{equation}\label{eq5}
J_{n}(2 x) = (\hat{c}\, x)^{n} e^{-\hat{c}\, x^{2}}\, \varphi(0), \quad \varphi(n) = \frac{1}{\Gamma(1+n)}.
\end{equation}

According to the above restyling, and by taking the freedom of treating $\hat{c}$ as an ordinary constant, an interesting plethora of new results concerning integrals of Bessel functions and their derivatives can be obtained. Just to give an example of how the method works, we consider the integral
\begin{equation}\label{eq6}
B_{\nu}(\alpha, \beta) = \int_{0}^{\infty} x\, J_{2\nu}(\alpha x) e^{i\beta x^{2}}\, dx, \quad (\alpha > 0, \, \beta > 0).
\end{equation}
The use of Eq.~(\ref{eq5}) allows to treat this integral as an elementary integral of exponential type and indeed, under the hypothesis $\alpha^{2} < 4 \beta$, we find
\begin{equation}\label{eq7}
B_{\nu}(\alpha, \beta) = \frac{1}{2}\, \left(\frac{\alpha}{2}\right)^{2\nu}\, \left(\frac{i}{\beta}\right)^{\nu + 1} b_{\nu}\left(-i\frac{\alpha^{2}}{4\beta}\right),
\end{equation}
where we have introduced the function
\begin{equation}\nonumber
b_{\nu}(x) = \sum_{k=0}^{\infty} \frac{\Gamma(\nu + k + 1)}{\Gamma(2\nu + k + 1)}\, \frac{x^{k}}{k!} = \frac{\sqrt{\pi}}{2} x^{\frac{1}{2}-\nu} e^{\frac{x}{2}}\, \left[I_{\nu-\frac{1}{2}}\left(\frac{x}{2}\right) + I_{\nu+\frac{1}{2}}\left(\frac{x}{2}\right)\right],
\end{equation}
where $I_{\mu}(z)$ is the modified Bessel function \cite{GEAndrews01}. In the case $\nu = 0$, Eq.~(\ref{eq7}) gives
\begin{equation}\label{eq8}
B_{0}(\alpha, \beta) = \frac{i}{2\beta} \exp\left(-i \frac{\alpha^{2}}{4\beta}\right),
\end{equation}
in agreement with Eqs.~\textbf{6.728.3} and \textbf{6.728.4} of \cite{ISGradsheteyn94}, and Eq.~\textbf{2.12.18.7} of \cite{APPrudnikov92}.

We consider now the Struve function \cite{GEAndrews01}
\begin{equation}\label{eq9}
\mathbf{H}_{\nu}(x) = \sum_{k=0}^{\infty} \frac{(-1)^{k}\, \left(\frac{x}{2}\right)^{2k+\nu+1}}{\Gamma\left(k + \frac{3}{2}\right) \Gamma\left(k + \nu + \frac{3}{2}\right)},
\end{equation}
that can be rewritten with the umbral operator $\hat{h}_{\nu}$ as
\begin{equation}\label{eq10}
\mathbf{H}_{\nu}(x) = \left(\frac{x}{2}\right)^{\nu+1}\, \exp\left[-\hat{h}_{\nu}\, \left(\frac{x}{2}\right)^{2}\right]\, \varphi(0),
\end{equation}
with
\begin{equation}\label{eq11}
\varphi(n) = (\hat{h}_{\nu})^{n} \varphi(0) = \frac{\Gamma(n+1)}{\Gamma\left(n + \frac{3}{2}\right)\, \Gamma\left(n + \nu + \frac{3}{2}\right)}.
\end{equation}
By applying our method, and using the Euler's reflection formula for Gamma functions, we find:
\begin{equation}\label{eq12}
\int_{0}^{\infty} \mathbf{H}_{\nu}(b x)\, dx = \frac{1}{b}\, \frac{\Gamma\left(1 + \frac{\nu}{2}\right)\, \Gamma\left(-\frac{\nu}{2}\right)}{\Gamma\left(\frac{1-\nu}{2}\right)\, \Gamma\left(\frac{1+\nu}{2}\right)} = - \frac{1}{b\, \tan\left(\frac{\pi\nu}{2}\right)}\, , \quad (-2 < \mathrm{Re}\nu < 0, b > 0),
\end{equation}
that is formula \textbf{6.811.1} of \cite{ISGradsheteyn94}. Further application of Eq.~(\ref{eq10}) furnishes
\begin{equation}\label{eq13}
\int_{-\infty}^{\infty} x^{-(\nu+1)}\, \mathbf{H}_{\nu}(x)\, dx = \frac{\sqrt{\pi}}{2^{\nu}}\, \hat{h}_{\nu}^{-\frac{1}{2}}\, \varphi(0) = \frac{\pi}{2^{\nu}\, \Gamma(1+\nu)}
\end{equation}
(see  formula \textbf{6.813.2} of \cite{ISGradsheteyn94}). It is evident that the method can easily be applied also to more complicated integrals involving the Struve functions.

As a further evidence of the usefulness of the method, we apply it to the search of close form of the following generating function
\begin{equation}\label{eq14}
G(x, t|\, m) = \sum_{n=0}^{\infty} \frac{t^{n}}{n!}\, J_{m\, n}(2 x).
\end{equation}
With the help of Eq.~(\ref{eq5}) it can be written as
\begin{equation}\label{eq15}
G(x, t|\, m) = \exp[(\hat{c} x)^{m}\, t - \hat{c}\, x^{2}]\, \varphi(0).
\end{equation}
If we introduce the higher, ($m \geq 2$), order Hermite polynomials \cite{FGTricomi54}
\begin{equation}\label{eq16}
H_{n}^{(m)}(u, v) = n! \sum_{k=0}^{[n/m]} \frac{u^{n-m k}\, v^{k}}{(n - m k)!\, k!},
\end{equation}
their generating function is \cite{FGTricomi54, PEAppell26, GDattoli97}:
\begin{equation}\label{eq17}
\sum_{n=0}^{\infty} \frac{z^{n}}{n!} H_{n}^{(m)}(u, v) = \exp(u z + v z^{m}).
\end{equation}
(Note that $H_{n}^{(2)}(x, y) = (-i)^{n}\, y^{n/2}\, H_{n}\left(\frac{i x}{2 \sqrt{y}}\right)$, where $H_{n}(z)$ are conventional Hermite polynomials \cite{GEAndrews01}). This allows us to recast Eq.~(\ref{eq15}) in the operator form
\begin{equation}\label{eq18}
G(x, t|\, m) = \sum_{n=0}^{\infty} \frac{\hat{c}^{\,n}}{n!} H_{n}^{(m)}(-x^{2}, x^{\,m}\, t) \varphi(0).
\end{equation}
Furthermore, as a consequence of the definitions (\ref{eq5}) and (\ref{eq16}), it is possible to show that
\begin{equation}\label{eq19}
G(x, t|\, m) = \,_{H}C_{0}^{\,(m)}\left[x^{2}, (-x)^{m} t\right],
\end{equation}
where
\begin{equation}\label{eq20}
_{H}C_{n}^{\,(m)}(x, y) = \sum_{k=0}^{\infty} \frac{(-1)^{k}}{k!\, (n + k)!}\, H_{k}^{(m)}(x, y),
\end{equation}
are the so called Hermite-based Tricomi functions \cite{GDattoli07}. The correctness of the above identity has been thoroughly checked numerically. We stress that its derivation with conventional means is much more involved.

In the same spirit, we will use operational methods to evaluate various families of integrals. In particular, we consider the integrals of the type
\begin{equation}\label{eq21}
I(x) = \int_{-\infty}^{\infty} f(x\, g(t))\, dt. 
\end{equation}
The identity \cite{GDattoli97}
\begin{equation}\label{eq22}
[g(t)]^{x\, \partial_{x}} f(x) = f\left[x\, g(t)\right]
\end{equation}
can be used to write
\begin{equation}\label{eq23}
I(x) = \int_{-\infty}^{\infty} [g(t)]^{x\, \partial_{x}} f(x)\, dt,
\end{equation}
and, by assuming that the integral ($a = \mathrm{constant}$)
\begin{equation}\nonumber
F(a) = \int_{-\infty}^{\infty} [g(t)]^{a}\, dt
\end{equation}
exists, one can express $I(x)$ as
\begin{equation}\label{eq24}
I(x) = \hat{F}(x\, \partial_{x}) f(x).
\end{equation}
If the function $f(x)$ admits a power series expansion
\begin{equation}\label{eq25}
f(x) = \sum_{k=0}^{\infty} \alpha(k) x^{m k + p},
\end{equation}
as a consequence of the identity
\begin{equation}\label{eq26}
\hat{F}(x\, \partial_{x}) x^{n} = F(n) x^{n},
\end{equation}
we get
\begin{equation}\label{eq27}
I(x) = \sum_{k=0}^{\infty} \alpha(k) F(m k + p) x^{m k + p}.
\end{equation}

To clarify the content of this result, we choose the case $g(t) = e^{t^{2}}$ , $f(x) = J_{n}(x)$. By using Eq.~(26) one has
\begin{equation}\label{eq28}
\int_{-\infty}^{\infty} J_{n}\left(x e^{-t^{2}}\right)\, dt = \sqrt{\pi}\, \sum_{k=0}^{\infty}\, \frac{(-1)^{k}}{k!\, (k + n)!}\, \frac{\left(\frac{x}{2}\right)^{2 k + n}}{\sqrt{2 k + n}}, \quad n > 0.
\end{equation}
The above series has been proved to be convergent and the correctness of the proposed procedure has been checked numerically. 

As a further example we consider the evaluation of the integral (\ref{eq21})
for $g(t) = (1 + t^{2})^{-1}$, $f(x) = x^{2} e^{-x^{2}}$. Since
\begin{equation}\label{eq29}
\int_{-\infty}^{\infty} (1 + t^{2})^{-a}\, dt = \sqrt{\pi}\, \frac{\Gamma\left(a - \frac{1}{2}\right)}{\Gamma(a)}
\end{equation}
we get
\begin{equation}\label{eq30}
\int_{-\infty}^{\infty} \frac{\exp\left[-\frac{x^{2}}{(1 + t^{2})^{2}}\right]}{(1 + t^{2})^{2}}\, dt = \sqrt{\pi}\, \sum_{k=0}^{\infty} \frac{(-x^{2})^{k}}{k!}\, \frac{\Gamma\left(2 k + \frac{3}{2}\right)}{2 k + 2} = \frac{\pi}{2} \,_{2}F_{2}\left(\frac{3}{4}, \frac{5}{4};\, 1, \frac{3}{2};\, -x^{2}\right),
\end{equation}
where $_{p}F_{q}(a_{1}, \ldots, a_{p};\, b_{1}, \ldots, b_{q};\, y)$ is the generalized hypergeometric function, see \cite{APPrudnikov92, GEAndrews01}.

Let us now consider the following Laplace-type transform integral
\begin{equation}\label{eq31}
L(x) = \int_{0}^{\infty} e^{-t}\, g(x\, t)\, dt\, ,
\end{equation}
also known as the Borel transform. The use of Eq.~(\ref{eq22}) leads to the operatorial identity
\begin{equation}\label{eq32}
L(x) = \Gamma(x\, \partial_{x} + 1)\, g(x)
\end{equation}
from which, under the assumption that $L(x)$ can be expanded as in Eq.~(\ref{eq1}), as a consequence of Eq.~(\ref{eq26}), one obtains
\begin{equation}\label{eq33}
g(x) = \sum_{k=0}^{\infty}\, \frac{\varphi(k)}{(k!)^{2}}\, (-x)^{k}.
\end{equation}
See Chap. 2 of \cite{MMDjrbashian93} for the use of the operator $\Gamma(x\, \partial_{x} + 1)$.

The use of Borel transform is of noticeable importance in the derivation of Quantum Chromodynamics (QCD) sum rules \cite{BLIoffe81}. For a general introduction see \cite{BLangwallner05}. We will just touch on this topic and leave further developments for a forthcoming paper. Just to quote a few examples we note that the polynomials given in Eq.~(\ref{eq16}) are the Borel transform of their hybrid counterpart \cite{GDattoli04}
\begin{equation}\label{eq34}
\tilde{H}_{n}^{(m)}(x, y) = \sum_{k=0}^{[n/m]}\, \frac{x^{n - m k}\, y^{k}}{k!\, [(n - m k)!]^{2}}, \quad m \geq 2,
\end{equation}
if we choose $g(x\, t) = \tilde{H}_{n}^{(m)}(x\, t, y)$ and treat $y$ as a parameter. They are called hybrid because they have properties in between those of Hermite and Laguerre polynomials. It is worth stressing that the Borel transform with respect to the $y$ variable [namely choosing $g(y t) = \tilde{H}_{n}^{(m)}(x, y t)\,$], is $e_{n}^{(m)}(x, y) = \sum_{k=0}^{[n/m]} \frac{x^{n - m k}\, y^{k}}{[(n - mk)!]^{2}}$ corresponding to a family  of truncated polynomials. See \cite{GDattoli03, GDattoli02} for discussions of related polynomial families. Furthermore, for $g(x) = (1 + x^{m})^{-1}$, $|x| < 1$, we find that the associated transform is the pseudo-trigonometric function \cite{PERicci78} $c_{0}^{(m)}(x)$, where
\begin{equation}\label{eq35}
c_{k}^{(m)}(x) = \sum_{r=0}^{\infty} (-1)^{r}\, \frac{x^{m r + k}}{(m r + k)!}, \quad 0 \leq k < m.
\end{equation}

The following relationship can be viewed as an extension of the Borel transform
\begin{equation}\label{eq36}
I_{\alpha, \beta}(x) = \int_{0}^{1} u^{\alpha - 1}\, (1 - u)^{\beta - 1}\, f(u\, x)\, du \quad (\mathrm{Re} \alpha, \mathrm{Re} \beta > 0).
\end{equation}
By using Eq.~(\ref{eq22}), this integral can be written in terms of the Beta function as follows
\begin{equation}\nonumber
I_{\alpha, \beta}(x) = B(\alpha + x\, \partial_{x}, \beta) f(x),
\end{equation}
i.e., in terms of the Gamma function
\begin{equation}\label{eq37}
I_{\alpha, \beta}(x) = \frac{\Gamma(\beta)\, \Gamma(\alpha + x\, \partial_{x})}{\Gamma(\alpha + \beta + x\, \partial_{x})} f(x).
\end{equation}
If the function $f(x)$ satisfies Eq.~(\ref{eq1}), one has
\begin{equation}\label{eq38}
I_{\alpha, \beta}(x) = \sum_{n=0}^{\infty} \frac{\Psi(n)}{n!} (-x)^{n}, 
\end{equation}
where
\begin{equation}\label{eq39}
\Psi(n) = B(\alpha + n, \beta)\, \varphi(n).
\end{equation}

In conclusion, in this paper we have provided a few elements proving the usefulness of formal procedures to get results encompassing various aspects of operators and of special functions. Other examples could be discussed to corroborate the validity of the method. We believe, however, that the main plot of the whole technique can be quite well understood via the examples discussed here. In closing, we remark that the results we have obtained are at the level of sound conjectures. The relevant rigorous proof of them is out of the scope of the present article. There are similar techniques that kept being revived, see for example \cite{IGonzalez10} and references therein.

\ \\

\noindent
\textbf{Acknowledgements}\\
One of us (G. D.) recognizes the warm hospitality and the financial support of the University of Paris XIII, whose stimulating atmosphere provided the necessary conditions for the elaboration of the ideas leading to this paper. K. G. thanks Funda\c{c}\~{a}o de Amparo \'{a} Pesquisa do Estado de S\~{a}o Paulo (FAPESP, Brazil) under Program No.~2010/15698-5. The authors acknowledge support from Agence Nationale de la Recherche (Paris, France) under Program PHYSCOMB No.~ANR-08-BLAN-243-2.

\end{document}